\documentclass[a4paper]{styles/svproc}
\usepackage{url}

\usepackage{amsmath} 
\usepackage{amssymb}  
\usepackage{graphicx}
\usepackage{algorithm}
\usepackage{algpseudocode}

\algrenewcommand\algorithmicensure{\textbf{Output:}}

\allowdisplaybreaks

\DeclareMathOperator*{\argmin}{argmin}

\DeclareMathOperator{\diag}{diag}

\begin{document}
\mainmatter              
\title{Square Root Gauss-Newton iLQR}
\titlerunning{Square Root Gauss-Newton iLQR}  
%
\author{Maximilian Haas-Heger \and Jur van den Berg}
\authorrunning{Maximilian Haas-Heger \and Jur van den Berg} 
%
\tocauthor{Maximilian Haas-Heger, Jur van den Berg}
\institute{Waabi Innovation Inc.}

\maketitle              

\begin{abstract}
The iterative Linear Quadratic Regulator (iLQR) is a widely used algorithm for nonlinear trajectory optimization. At each iteration, it solves a local linear-quadratic approximation of the problem via dynamic programming, propagating a quadratic cost-to-go function. If the Hessian of the cost-to-go approximation is positive-semidefinite, one can derive a \emph{square root} formulation of iLQR that propagates its Cholesky factor instead. This offers significant numerical advantages\textemdash much as square root Kalman filters improve upon their conventional counterparts\textemdash particularly when iLQR is used within an augmented Lagrangian framework for handling constraints, where large penalties degrade conditioning. Previous square root formulations of iLQR and related algorithms exist, but they are either numerically suboptimal, algorithmically complex, or both. In this paper, we show that the key to an effective square root formulation lies in the Gauss-Newton (weighted least-squares) structure of the cost function: this yields a positive semidefiniteness property that extends beyond the Hessian to the full augmented cost-to-go matrix, and enables a backward pass of remarkable simplicity in which each step reduces to a single QR-decomposition, from which the feedback gain and propagated Cholesky factor are extracted directly.
\end{abstract}

\section{Introduction}
Trajectory optimization\textemdash the problem of finding a sequence of controls that steers a dynamical system along a low-cost path\textemdash is a cornerstone of robotics, autonomous driving, and aerospace planning. When the dynamics are linear and the cost is quadratic with positive-(semi)definite Hessians, the problem can be solved exactly with a Linear Quadratic Regulator (LQR). In general, however, the dynamics are nonlinear and the cost is non-quadratic, and we must resort to iterative methods. The most general such approach is differential dynamic programming (DDP)~\cite{jacobsen}, in which a full second-order expansion of the action-value function is performed in each step. 
The \emph{iterative Linear Quadratic Regulator} (iLQR)~\cite{li,todorov} closely resembles DDP but applies a Gauss-Newton simplification by linearizing the dynamics instead  \cite{roulet,sideris}. 

The key idea behind iLQR is conceptually simple: given a current candidate trajectory, linearize the dynamics and quadratize the cost around that trajectory, solve the resulting linear-quadratic subproblem exactly, and then use the solution to obtain an improved trajectory. This process repeats until convergence. The linear-quadratic subproblem itself is solved efficiently via \emph{dynamic programming}: starting from the end of the trajectory and working backward, the algorithm computes a quadratic \emph{cost-to-go} (or value) function for each time step, along with a locally optimal linear feedback policy.

In standard LQR, the Hessian of the value function is \emph{positive-semidefinite} (like the covariance matrix is in a Kalman filter). This property applies to iLQR as well if the cost is expressed as a weighted sum of squared \emph{residuals}. The Gauss-Newton nature of iLQR can then be extended to the cost function by linearizing the residuals rather than quadratizing the cost as a whole \cite{geoffroy,giftthaler,sideris}. This not only further simplifies the algorithm, but also enables the formulation of a \emph{square root} version of iLQR that\textemdash analogous to a square root Kalman filter~\cite{bierman,chin}\textemdash propagates the Cholesky factor of the cost-to-go Hessian rather than the matrix itself. This offers significant numerical advantages, particularly in poorly conditioned problems and where the available precision is limited, as representing the cost-to-go via its Cholesky factor roughly doubles the effective precision~\cite{dewilde} and guarantees numerical positive-semidefiniteness by construction.

The square root formulation becomes especially valuable when iLQR is combined with an \emph{augmented Lagrangian} framework for handling constraints
. In this setting, constraint violations are penalized by squared residual terms weighted by a large penalty parameter $\mu$ 
that widen the eigenvalue spread in the cost-to-go matrices and degrade the conditioning of iLQR. The square root formulation mitigates this directly, as it operates on Cholesky factors whose condition numbers scale with $\sqrt{\mu}$ rather than $\mu$.

Whereas square root versions of the Kalman filter have enjoyed widespread adoption (including its storied implementation aboard the Apollo spacecraft) \cite{chin,grewal}, they have gained relatively limited traction within the context of iLQR. This is somewhat surprising, as linear-quadratic estimation and control are in many ways each other's dual \cite{todorov08}. The work of \cite{geoffroy} suggests a square root version of DDP, also with sum-of-squares cost functions, that propagates a non-triangular square root of the value function whose dimension grows with each time step, and from which the control policy is extracted using a pseudo-inverse. Square root techniques have also been proposed for solvers of algebraic Riccati equations \cite{dewilde,frison13,frison14}, which are closely related to linear-quadratic control, and for the Riccati differential equation \cite{tedrake}, where numerical integration can cause loss of symmetry or positive-definiteness. Perhaps the most direct attempt at creating a (triangular) square root formulation of iLQR is that of~\cite{altro}, which also applies iLQR within an augmented Lagrangian framework. 
The approach is suboptimal, however, both numerically and in terms of complexity, as it requires three separate QR-decompositions and a series of rank-one ``downdates'' per time step.

In this paper, we present a square root iLQR approach that avoids these shortcomings. We exploit an underappreciated benefit of the Gauss-Newton approximation in that not only the Hessian but the full augmented matrix of the cost-to-go function is positive-semidefinite. This enables a square root backward pass of remarkable simplicity and elegance: each step reduces to performing a single QR-decomposition, from which the feedback gain and propagated Cholesky factor are extracted directly. We will demonstrate the effectiveness and robustness of our square root iLQR implementation using numerical experiments.

\section{Problem Formulation}\label{sec:problem}
Consider a discrete-time dynamical system whose evolution is governed by:
\begin{align}
\mathbf{x}_{t + 1} & = \mathbf{f}_t(\mathbf{x}_t, \mathbf{u}_t),
\end{align}
where $\mathbf{x}_t \in \mathbb{R}^n$ is the \emph{state} (e.g., position, velocity, orientation of a robot) and $\mathbf{u}_t \in \mathbb{R}^m$ is the \emph{control input} (e.g., acceleration, steering) at time $t$. The function $\mathbf{f}_t$ encodes how the system transitions from one state to the next given a control.

Our goal is to choose a sequence of controls $\mathbf{u}_0, \ldots, \mathbf{u}_{\ell - 1}$ over a horizon of $\ell$ steps that minimizes the total cost:
\begin{align}
\text{minimize} &\quad c_\ell(\mathbf{x}_\ell) + \sum_{t = 0}^{\ell - 1} c_t(\mathbf{x}_t, \mathbf{u}_t) \label{eq:totalcost} \\
\text{subject to} &\quad \mathbf{x}_{t+1} = \mathbf{f}_t(\mathbf{x}_t, \mathbf{u}_t), \quad \mathbf{x}_0 = \mathbf{x}_0^\star,   \nonumber
\end{align}
where the states $\mathbf{x}_0,\ldots,\mathbf{x}_\ell$ are determined by the given initial state $\mathbf{x}^\star_0$ and the chosen controls through the dynamics function $\mathbf{f}_t$. Here, $c_t(\mathbf{x}_t, \mathbf{u}_t)$ is the \emph{running cost} incurred at time $t$ (e.g., penalizing deviation from a desired path or excessive control effort), and $c_\ell(\mathbf{x}_\ell)$ is the \emph{terminal cost} that applies to the final state. In this paper, we focus on cost functions of the \emph{weighted least-squares} form:
\begin{align}
c_t(\mathbf{x}_t, \mathbf{u}_t) & = \frac{1}{2}\mathbf{r}_t(\mathbf{x}_t, \mathbf{u}_t)^T W_t \mathbf{r}_t(\mathbf{x}_t, \mathbf{u}_t), & W_t & \succ 0, \label{eq:gaussnewtoncost} \\
c_\ell(\mathbf{x}_\ell) & = \frac{1}{2}\mathbf{r}_\ell(\mathbf{x}_\ell)^T W_\ell \mathbf{r}_\ell(\mathbf{x}_\ell), & W_\ell &\succ 0. \label{eq:gaussnewtonterminal}
\end{align}
where $W_t \in \mathbb{R}^{k\times k}$ is a positive-definite weight matrix and $\mathbf{r}_t \in \mathbb{R}^k$ is a vector of \emph{residuals}. In practice, many trajectory optimization objectives can be expressed in this form \cite{geoffroy,giftthaler}.

A problem of this type can generally be expressed in terms of the Bellman equation, which recursively defines a value function (or cost-to-go function) $v_t(\mathbf{x}_t)$ that represents the minimum future cost achievable starting from state $\mathbf{x}_t$ at time $t$. For the final stage $\ell$, we have:
\begin{align}
v_\ell(\mathbf{x}_\ell) & = c_\ell(\mathbf{x}_\ell), \label{eq:vell}
\end{align}
and then backward in time for $t = \ell-1,\ldots,0$:
\begin{align}
q_t(\mathbf{x}_t, \mathbf{u}_t) & = c_t(\mathbf{x}_t, \mathbf{u}_t) + v_{t + 1}(\mathbf{f}_t(\mathbf{x}_t, \mathbf{u}_t)), \label{eq:qfunction} \\
\pi_t(\mathbf{x}_t) & = \argmin_{\mathbf{u}_t} q_t(\mathbf{x}_t, \mathbf{u}_t), \label{eq:policy} \\
v_t(\mathbf{x}_t) & = \min_{\mathbf{u}_t} q_t(\mathbf{x}_t, \mathbf{u}_t) = q_t(\mathbf{x}_t, \pi_t(\mathbf{x}_t)). \label{eq:vt}
\end{align}
The intermediate functions $q_t$ and $\pi_t$ have natural interpretations: $q_t(\mathbf{x}_t, \mathbf{u}_t)$ is the \emph{action-value} function, which represents the cost of taking action $\mathbf{u}_t$ in state $\mathbf{x}_t$ and then acting optimally thereafter. The function $\pi_t(\mathbf{x}_t)$ is the \emph{policy}, which gives the optimal action to take at time $t$ given the state $\mathbf{x}_t$.

\section{Standard Gauss-Newton iLQR}\label{sec:standard}

The iLQR algorithm assumes that we are given a current \emph{nominal trajectory} consisting of states $\mathbf{x}^\star_0,\ldots,\mathbf{x}^\star_\ell$ and controls $\mathbf{u}^\star_0,\ldots,\mathbf{u}^\star_{\ell - 1}$ that satisfies the dynamics constraint, i.e., $\mathbf{x}^\star_{t+1} = \mathbf{f}_t(\mathbf{x}^\star_t, \mathbf{u}^\star_t)$, for all $t$. We work in terms of small deviations from this trajectory:
\begin{align}
\delta \mathbf{x}_t &= \mathbf{x}_t - \mathbf{x}^\star_t, & \delta \mathbf{u}_t &= \mathbf{u}_t - \mathbf{u}^\star_t.
\end{align}

The algorithm repeatedly performs two passes over the trajectory:
\begin{enumerate}
    \item \textbf{Backward pass.} Starting from the terminal state and moving backward in time, locally approximate the dynamics and each residual vector as linear around the current trajectory, and solve the resulting linear-quadratic subproblem via dynamic programming. This produces a linear feedback policy at each step.
    \item \textbf{Forward pass.} Starting from the initial state, simulate the system forward using the linear feedback policy to obtain a new, improved trajectory.
\end{enumerate}
These two passes alternate until the trajectory converges. 

\subsection{Local approximations}
In iLQR the dynamics and cost are approximated as linear and quadratic functions, respectively, around the current nominal trajectory.

\textbf{Dynamics linearization:} We linearize the dynamics around the nominal trajectory as
$\mathbf{f}_t(\mathbf{x}_t, \mathbf{u}_t) \approx \mathbf{f}_t(\mathbf{x}^\star_t, \mathbf{u}^\star _t) + \nabla_\mathbf{x}\mathbf{f}_t\, \delta\mathbf{x}_t + \nabla_\mathbf{u}\mathbf{f}_t\, \delta\mathbf{u}_t$.
Recognizing that $\mathbf{f}_t(\mathbf{x}_t, \mathbf{u}_t) - \mathbf{f}_t(\mathbf{x}^\star_t, \mathbf{u}^\star _t) = \mathbf{x}_{t + 1} - \mathbf{x}_{t + 1}^\star = \delta \mathbf{x}_{t + 1}$, we can write this succinctly as:
\begin{align}
\begin{bmatrix} \delta \mathbf{x}_{t + 1} \\ 1 \end{bmatrix} & \approx F_t \begin{bmatrix} \delta \mathbf{u}_t \\ \delta \mathbf{x}_t \\ 1 \end{bmatrix}, & F_t & = \left. \begin{bmatrix} \nabla_\mathbf{u}\mathbf{f}_t & \nabla_\mathbf{x}\mathbf{f}_t & \mathbf{0} \\ \mathbf{0}^T & \mathbf{0}^T & 1 \end{bmatrix}\right\rvert_{\mathbf{x}^\star_t, \mathbf{u}^\star_t}, \label{eq:lindynamics}
\end{align}
where $F_t$ is a matrix of Jacobians evaluated at the nominal trajectory. 

\textbf{Cost quadratization:} The running cost is approximated as a quadratic function around the nominal trajectory:\footnote{Appending a ``1'' to the state vector lets a quadratic with linear and constant terms be written as a single matrix expression. We call the matrix involved an \emph{augmented matrix}, and denote such matrices by majuscule bold symbols throughout.}
\begin{align}
c_t(\mathbf{x}_t, \mathbf{u}_t) & \approx \frac{1}{2}\begin{bmatrix} \delta\mathbf{u}_t \\ \delta\mathbf{x}_t \\ 1 \end{bmatrix}^T \mathbf{C}_t \begin{bmatrix} \delta\mathbf{u}_t \\ \delta\mathbf{x}_t \\ 1 \end{bmatrix}, & 
c_\ell(\mathbf{x}_\ell) & \approx \frac{1}{2}\begin{bmatrix} \delta\mathbf{x}_\ell \\ 1 \end{bmatrix}^T \mathbf{C}_\ell \begin{bmatrix} \delta\mathbf{x}_\ell \\ 1 \end{bmatrix}.  \label{eq:quadcost}
\end{align}
Rather than quadratizing the cost directly, i.e. expanding
$c_t$ to second order, the Gauss-Newton version of iLQR approximates the residual vector $\mathbf{r}_t$ linearly:
\begin{align}
\mathbf{r}_t(\mathbf{x}_t, \mathbf{u}_t) & \approx R_t \begin{bmatrix}\delta\mathbf{u}_t \\\delta\mathbf{x}_t \\ 1 \end{bmatrix}, & \mathbf{r}_\ell(\mathbf{x}_\ell ) & \approx R_\ell \begin{bmatrix}\delta\mathbf{x}_\ell \\ 1 \end{bmatrix}, \label{eq:linresidual}
\end{align}
with Jacobians and residual values evaluated at the nominal trajectory:
\begin{align}
R_t & = \left. \begin{bmatrix} \nabla_\mathbf{u}\mathbf{r}_t & \nabla_\mathbf{x}\mathbf{r}_t & \mathbf{r}_t \end{bmatrix} \right\rvert_{\mathbf{x}^\star_t,\mathbf{u}^\star_t}, & R_\ell & = \left. \begin{bmatrix} \nabla_\mathbf{x}\mathbf{r}_\ell & \mathbf{r}_\ell \end{bmatrix} \right\rvert_{\mathbf{x}^\star_\ell,\mathbf{u}^\star_\ell}.
\end{align}
Substituting Eq.\ \eqref{eq:linresidual} into Eq.\ \eqref{eq:gaussnewtoncost} and comparing with Eq.\ \eqref{eq:quadcost} gives:
\begin{align}
\mathbf{C}_t & = R_t^T W_t R_t \succeq 0, \label{eq:cmatrix}
\end{align}
which is positive-semidefinite by construction.

\subsection{Backward pass}
Given linear dynamics and quadratic cost, the cost-to-go function $v_t$ as defined in Eqs.\ \eqref{eq:vell}--\eqref{eq:vt} and policy $\pi_t$ as defined in Eq.\ \eqref{eq:policy} are quadratic and linear, respectively, for all $t$, and can be explicitly evaluated in closed-form. Let the (approximate) value function be defined relative to the nominal trajectory as:
\begin{align}
v_t(\mathbf{x}_t) & \approx \frac{1}{2}\begin{bmatrix} \delta\mathbf{x}_t \\ 1 \end{bmatrix}^T \mathbf{S}_t \begin{bmatrix} \delta \mathbf{x}_t \\ 1 \end{bmatrix}, & \mathbf{S}_t & = \begin{bmatrix} S_\mathbf{xx} & \mathbf{s}_\mathbf{x} \\ \mathbf{s}_\mathbf{x}^T & s \end{bmatrix}, \label{eq:costogo}
\end{align}
where $S_\mathbf{xx}$ approximates the Hessian of the cost-to-go function, $\mathbf{s}_\mathbf{x}$ its gradient, and $s$ encodes twice its value at the nominal trajectory. Then initially, for $t = \ell$, since the terminal cost is approximated by a quadratic function as well (see Eq.\ \eqref{eq:quadcost} and Eq.\ \eqref{eq:cmatrix}), we have, following Eq.\ \eqref{eq:vell}:
\begin{align}
\mathbf{S}_\ell & = R^T_\ell W_\ell R_\ell \succeq 0. \label{eq:sell}
\end{align}
From this starting point, the algorithm steps backward in time, computing $\pi_t$ (as defined in Eq.\ \eqref{eq:policy}) and $v_t$ for each $t = \ell - 1, \ldots, 0$.

Combining the cost-to-go function of time $t + 1$ with the local approximations of time $t$, we obtain the \emph{action-value} function at time $t$ as defined in Eq.~\eqref{eq:qfunction}:
\begin{align}
q_t(\mathbf{x}_t, \mathbf{u}_t) &= \frac{1}{2}\begin{bmatrix} \delta \mathbf{u}_t \\ \delta \mathbf{x}_t \\ 1 \end{bmatrix}^T \mathbf{Q}_t \begin{bmatrix} \delta \mathbf{u}_t \\ \delta \mathbf{x}_t \\ 1 \end{bmatrix}, & \mathbf{Q}_t & = \begin{bmatrix} Q_\mathbf{uu} & Q_\mathbf{ux} & \mathbf{q}_\mathbf{u} \\ Q_\mathbf{ux}^T & Q_\mathbf{xx} & \mathbf{q}_\mathbf{x} \\ \mathbf{q}_\mathbf{u}^T & \mathbf{q}_\mathbf{x}^T & q \end{bmatrix},
\label{eq:qmatrix}
\end{align}
where the matrix $\mathbf{Q}_t$ is obtained by substituting the linearized dynamics of Eq.\ \eqref{eq:lindynamics} into the cost-to-go function of time $t+1$ (Eq.\ \eqref{eq:costogo}) and adding the immediate cost of the current time step (Eq.\ \eqref{eq:quadcost} and Eq.\ \eqref{eq:cmatrix}):
\begin{align}
\mathbf{Q}_t & = R^T_t W_t R_t + F_t^T \mathbf{S}_{t + 1} F_t \succeq 0. \label{eq:qt}
\end{align}

To find the optimal control, we minimize $q_t$ over $\delta\mathbf{u}_t$ (see Eq.\ \eqref{eq:policy}). Since $q_t$ is quadratic, we can do so in closed form if $Q_\mathbf{uu} \succ 0$; Section~\ref{sec:advantages} discusses what guarantees this. 
 Solving $\partial q_t / \partial \delta\mathbf{u}_t = 0$ gives:
\begin{align}
    \delta\mathbf{u}_t & = -Q_{\mathbf{uu}}^{-1}\begin{bmatrix} Q_\mathbf{ux} & \mathbf{q}_\mathbf{u} \end{bmatrix}\begin{bmatrix} \delta \mathbf{x}_t \\ 1 \end{bmatrix} = -L_t \begin{bmatrix} \delta \mathbf{x}_t \\ 1 \end{bmatrix}, & L_t & = Q_\mathbf{uu}^{-1} \begin{bmatrix} Q_\mathbf{ux} & \mathbf{q}_\mathbf{u} \end{bmatrix}, \label{eq:matrixlt}
\end{align}
where
$L_t$
is the \emph{feedback gain} matrix defining the (locally) optimal policy $\pi_t$. 

Substituting the optimal $\delta\mathbf{u}_t$ back into Eq.~\eqref{eq:qmatrix} eliminates the control and gives the updated cost-to-go function at time $t$ (see Eqs.\ \eqref{eq:vt} and \eqref{eq:costogo}) with:
\begin{align}
\mathbf{S}_t & = \begin{bmatrix} Q_\mathbf{xx} & \mathbf{q}_\mathbf{x} \\ \mathbf{q}_\mathbf{x}^T & q \end{bmatrix} - L_t^T Q_\mathbf{uu} L_t \succeq 0. \label{eq:matrixst}
\end{align}
This is the well-known \emph{Riccati recursion}. The subtraction of $L_t^T Q_\mathbf{uu} L_t \succeq 0$ represents the cost reduction achieved relative to the nominal trajectory by choosing controls optimally. The recursion continues until $\mathbf{S}_0$ has been computed.

\subsection{Forward pass} \label{sec:forwardpass}
Once the backward pass is complete, we have a feedback policy $L_t$ for all $t$. The forward pass simulates the system from the given initial state $\mathbf{x}_0 \gets \mathbf{x}_0^\star$ using:
\begin{align}
\mathbf{u}_t &\gets \mathbf{u}^\star_t - L_t \begin{bmatrix} \mathbf{x}_t - \mathbf{x}_t^\star \\ \alpha \end{bmatrix}, & \mathbf{x}_{t+1} &\gets \mathbf{f}_t(\mathbf{x}_t, \mathbf{u}_t),
\end{align}
where $\alpha \in (0, 1]$ is a step-size parameter. When $\alpha = 1$, the full policy of the linear-quadratic subproblem is applied, while $\alpha = 0$ yields the current trajectory; any value in between effectively interpolates between these two extremes.

Let $j(\alpha)$ be the total cost of the trajectory resulting from the choice of $\alpha$, and let $j^\star = j(0)$ be the total cost of the current nominal trajectory. A \emph{backtracking line search} selects $\alpha$ by starting from $\alpha = 1$ and halving it until the Armijo-Goldstein sufficient-decrease condition is satisfied \cite{nocedal}:
\begin{align}
j(\alpha) - j^\star & \leq  \beta\,\alpha\,d, \label{eq:armijo}
\end{align}
where $\beta \in (0, \frac{1}{2}]$ and $d = dj/d\alpha|_{\alpha=0}$. The value $d < 0$ equals twice the cost change predicted by the linear-quadratic model,\footnote{This follows because the change in cost under the linear-quadratic model is a quadratic function of $\alpha$ that attains its minimum at $\alpha = 1$.} and is obtained in the backward pass by taking the sum over all $t$
of the bottom-right entry of the matrix $-L_t^T Q_\mathbf{uu} L_t$ (see Eq.\ \eqref{eq:matrixst}).\footnote{Were the constant term of the quadratized cost\textemdash the bottom-right entry of $\mathbf{C}_t$\textemdash set to zero, this same value would accumulate in the entry $s$ of $\mathbf{S}_t$ and could be extracted from $\mathbf{S}_0$.} The resulting trajectory becomes the new nominal, and the alternation of backward and forward passes continues until $|\frac{1}{2}d/j^\star|$ falls below a specified tolerance, indicating that the optimization has converged.

\subsection{Advantages of the Gauss-Newton Approach}\label{sec:advantages}
The Gauss-Newton specialization carries several important advantages over working with general cost functions:
\begin{enumerate}
\item \textbf{No second derivatives.} The approach requires only first derivatives (Jacobians) of the residual functions, not second derivatives of the cost.
\item \textbf{Positive-semidefiniteness by construction.} Since $\mathbf{C}_t = R_t^T W_t R_t$ with $W_t \succ 0$, the augmented matrix $\mathbf{C}_t$ is positive-semidefinite automatically. No regularization is needed. Moreover, the elements of $\mathbf{C}_t$ encoding the gradients and value are exact.
\item \textbf{The full augmented cost-to-go matrix is positive-semidefinite.} Since $\mathbf{C}_t \succeq 0$ and $\mathbf{S}_{t+1} \succeq 0$ implies $F_t^T \mathbf{S}_{t+1} F_t \succeq 0$, their sum $\mathbf{Q}_t = \mathbf{C}_t + F_t^T \mathbf{S}_{t+1} F_t \succeq 0$ (Eq.~\eqref{eq:qt}). In turn, $\mathbf{S}_t$ is the Schur complement of $Q_\mathbf{uu} \succ 0$ in $\mathbf{Q}_t \succeq 0$, so $\mathbf{S}_t \succeq 0$ as well. Since $\mathbf{S}_\ell = \mathbf{C}_\ell \succeq 0$, the argument applies by induction to all $t$. This property\textemdash positive-semidefiniteness of the \emph{full augmented} cost-to-go matrix, not just its Hessian block\textemdash is the key enabler of the square root formulation that follows in Section~\ref{sec:sqrt}.
\end{enumerate}

The one prerequisite this structure does not supply is $Q_\mathbf{uu} \succ 0$. Here $C_\mathbf{uu} = \nabla_\mathbf{u}^T\mathbf{r}_t\, W_t \nabla_\mathbf{u}\mathbf{r}_t$, which is positive-definite exactly when the Jacobian $\nabla_\mathbf{u}\mathbf{r}_t \in \mathbb{R}^{k \times m}$ has full column rank\textemdash possible only if $k \geq m$\textemdash and we impose this as a requirement on the residual functions. It is a mild one in practice, as a strictly positive cost on every control channel suffices.

Let us take a closer look at how the Gauss-Newton approach creates a positive-semidefinite approximation of the true Hessians of the running cost functions.
For a cost of the form $c(\mathbf{z}) = \frac{1}{2}\mathbf{r}(\mathbf{z})^T W \mathbf{r}(\mathbf{z})$, the true Hessian is:
\begin{align}
\nabla_{\mathbf{zz}}^2 c & = \nabla_\mathbf{z}^T\mathbf{r}\, W\, \nabla_\mathbf{z}\mathbf{r} + 
(W\mathbf{r}) \cdot \nabla^2_{\mathbf{zz}}\mathbf{r}.
 \label{eq:secondorderc}
\end{align}
The Gauss-Newton approximation retains only the first term. This is why it is particularly effective when either the residuals are small (making $W\mathbf{r} \approx \mathbf{0}$) or the residual functions are nearly linear (making $\nabla^2_{\mathbf{zz}}\mathbf{r} \approx 0$). This is also what sets the formulation apart from the usual ones. With general cost functions $Q_\mathbf{uu}$ is positive-definite only near an optimum, and iLQR variants force positivity by Levenberg-Marquardt regularization of $Q_\mathbf{uu}$ \cite{li,altro,todorov} or of $S_\mathbf{xx}$ \cite{ilqr}, by setting $Q_\mathbf{uu} = I$ \cite{roulet}, or by clamping negative eigenvalues to obtain the nearest positive-semidefinite Hessian \cite{higham88,eigenvalues}. Regularizing $Q_\mathbf{uu}$ alone, however, leaves the cost-to-go Hessians $S_\mathbf{xx}$ indefinite in general, and that precludes a square root formulation outright: an indefinite matrix has no real Cholesky factor.

\section{Square Root iLQR}\label{sec:sqrt}
The Gauss-Newton cost structure guarantees that the augmented matrices $\mathbf{S}_t$ and $\mathbf{Q}_t$ are positive-semidefinite for all $t$. This is what makes a square root formulation possible. Instead of propagating $\mathbf{S}_t$, the algorithm propagates an upper-triangular Cholesky factor $\sqrt{\mathbf{S}_t}$ satisfying $\sqrt{\mathbf{S}_t}\vphantom{Q}^{\,T}\!\sqrt{\mathbf{S}_t} = \mathbf{S}_t$, whose condition number equals the \emph{square root} of that of $\mathbf{S}_t$, roughly doubling the number of reliable significant digits.

In the backward pass of square root iLQR we wish to compute $\sqrt{\mathbf{S}_t}$ and $L_t$ given $\sqrt{\mathbf{S}_{t + 1}}$.
The derivation proceeds in two steps. First, we show that $\sqrt{\mathbf{S}_t}$ appears as the bottom-right block of the Cholesky factor $\sqrt{\mathbf{Q}_t}$ of $\mathbf{Q}_t$. Then, we show how to compute $\sqrt{\mathbf{Q}_t}$ directly from available inputs using a QR-decomposition, without ever explicitly constructing $\mathbf{Q}_t$ or any other squared matrix.


\subsection{Cholesky factorization of $\mathbf{Q}_t$}\label{sec:cholesky_schur}
Any positive-definite matrix $Q \succ 0$ has a \emph{unique} (up to sign conventions on the diagonal) upper-triangular Cholesky factor $\sqrt{Q}$ such that $Q = \sqrt{Q}\vphantom{Q}^{\,T}\!\sqrt{Q}$. If $Q$ is partitioned into two-by-two blocks, its Cholesky factor  is given by:
\begin{align}
Q & = \begin{bmatrix} Q_{00} & Q_{01} \\ Q_{01}^T & Q_{11} \end{bmatrix}, &
\sqrt{Q} & = \begin{bmatrix} \sqrt{Q_{00}} & \sqrt{Q_{00}}^{-T} Q_{01} \\ 0 & \sqrt{Q_{11} - Q_{01}^T Q_{00}^{-1} Q_{01}} \end{bmatrix}. \label{eq:cholesky}
\end{align}
The bottom-right entry involves $Q_{11} - Q_{01}^T Q_{00}^{-1} Q_{01}$, the \emph{Schur complement} of $Q_{00}$ in $Q$. For the matrix $\mathbf{Q}_t$ as defined in Eq.~\eqref{eq:qmatrix}, partitioned with $Q_{00}$ conformal to $Q_{\mathbf{uu}}$, this Schur complement is precisely $\mathbf{S}_t$ as given by Eq.\ \eqref{eq:matrixst}. If we let $\sqrt{\mathbf{Q}_t}  = \bigl[\begin{smallmatrix} U_{00} & U_{01} \\ 0 & U_{11} \end{smallmatrix}\bigr]$, then by Eq.\ \eqref{eq:matrixlt} and Eq.\ \eqref{eq:matrixst} the blocks are:
\begin{align}
U_{00} & = \sqrt{Q_{00}} = \sqrt{Q_{\mathbf{uu}}}, \\
U_{01} & = \sqrt{Q_{00}}^{-T} Q_{01} 
= \sqrt{Q_{00}}\, Q_{00}^{-1} Q_{01} = \sqrt{Q_{\mathbf{uu}}}\, L_t, \label{eq:u01} \\
U_{11} & = \sqrt{Q_{11} - Q_{01}^T Q_{00}^{-1} Q_{01}} 
       = \sqrt{\mathbf{S}_t}. \label{eq:u11}
\end{align}
Thus, the Cholesky factorization of $\mathbf{Q}_t$ becomes:
\begin{align}
\mathbf{Q}_t & = \sqrt{\mathbf{Q}_t}^T\!\sqrt{\mathbf{Q}_t}, & \sqrt{\mathbf{Q}_t} & = \begin{bmatrix} \sqrt{Q_\mathbf{uu}} & \sqrt{Q_\mathbf{uu}}\, L_t \\
0 & \sqrt{\mathbf{S}_t} \end{bmatrix}. \label{eq:sqrtq}
\end{align}

In other words, the quantities we seek\textemdash $\sqrt{\mathbf{S}_t}$ and $L_t$\textemdash are encoded directly in the Cholesky factor $\sqrt{\mathbf{Q}_t}$ of $\mathbf{Q}_t$. The remaining question is how to compute this Cholesky factor from square root inputs without forming any squared matrix.

\subsection{Computing $\sqrt{\mathbf{Q}_t}$ via QR-decomposition}\label{sec:computing_sqrt_q}
The second key fact underpinning square root iLQR is that if we can construct an alternative factorization $\mathbf{Q}_t = M^T M$ of $\mathbf{Q}_t$ for some matrix $M$, then the Cholesky factor $\sqrt{\mathbf{Q}_t}$ can be obtained directly from the QR-decomposition of $M$.
The QR-decomposition of a matrix $M \in \mathbb{R}^{a\times b}$ ($a \geq b$) produces a semi-orthogonal matrix $Q \in \mathbb{R}^{a\times b}$ ($Q^T Q = I$) and an upper-triangular matrix $U \in \mathbb{R}^{b\times b}$ such that $M = QU$ \cite{golub}. Since $Q$ is semi-orthogonal, we have $M^T M = U^T Q^T Q U = U^T U$. Thus, if $\mathbf{Q}_t = M^T M$, then by definition $U$ is the Cholesky factor $\sqrt{\mathbf{Q}_t}$, as it is upper-triangular and satisfies $U^T U = \mathbf{Q}_t$.

Such an alternative factorization of $\mathbf{Q}_t$ can in fact be constructed straightforwardly. Substituting Eq.~\eqref{eq:cmatrix} into Eq.~\eqref{eq:qt} gives $\mathbf{Q}_t = R_t^T W_t R_t + F_t^T \mathbf{S}_{t+1} F_t$, which can be rewritten as:
\begin{align}
\mathbf{Q}_t & = M_t^T M_t, & M_t & = \begin{bmatrix}
\sqrt{W_t}\, R_t \\
\sqrt{\mathbf{S}_{t + 1}}\, F_t
\end{bmatrix}. \label{eq:matrixm}
\end{align}

Hence, the QR-decomposition of $M_t$ yields $\sqrt{\mathbf{Q}_t}$ directly, without ever explicitly constructing $\mathbf{Q}_t$. Since QR-decomposition operates via orthogonal transformations (typically Householder reflections) that preserve norms and do not amplify rounding errors, computing $\sqrt{\mathbf{Q}_t}$ this way is numerically superior to explicitly computing $M_t^T M_t$ and then factoring it using a Cholesky decomposition.

Moreover, QR-factorization can robustly produce Cholesky factors of posi-tive-semidefinite matrices (such as $\mathbf{Q}_t$). The Cholesky factor is not unique in that case, but this does not affect the correctness of square root iLQR: Consider $\sqrt{\mathbf{Q}_t}$ as partitioned in Eq.\ \eqref{eq:sqrtq}. Since $Q_\mathbf{uu} \succ 0$, the blocks in the top row will be unique. The bottom-right block $\sqrt{\mathbf{S}_t}$ may not be unique as $\mathbf{S}_t \succeq 0$, but it will be a faithful Cholesky factor of $\mathbf{S}_t$ in the sense that it is upper-triangular and that $\sqrt{\mathbf{S}_t}\vphantom{Q}^{\,T}\!\sqrt{\mathbf{S}_t} = \mathbf{S}_t$.

Note, lastly, that we do not use the $Q$-factor of the QR-factorization. This allows for a particularly efficient ``Q-less'' implementation
that requires (for $a \geq b$) only $\sim 2ab^2 - \tfrac{2}{3}b^3$ floating-point operations (omitting sub-cubic terms) \cite{golub}.

\subsection{Summary}
Since $M_t$ is constructed from available inputs only (including $\sqrt{\mathbf{S}_{t + 1}}$), and $\sqrt{\mathbf{Q}_t}$ contains the quantities we seek ($\sqrt{\mathbf{S}_t}$ and $L_t$), each step of the backward pass in square root iLQR reduces to executing the following three simple steps:
\begin{enumerate}
\item Construct matrix $M_t$ as defined in Eq.\ \eqref{eq:matrixm} using $\sqrt{\mathbf{S}_{t+1}}$ (and other inputs).
\item \label{step:qr} Compute the upper-triangular matrix
\begin{align} U & = \sqrt{\mathbf{Q}_t} = \begin{bmatrix} U_{00} & U_{01} \\ 0 & U_{11} \end{bmatrix} \end{align}
through the (Q-less) QR-decomposition of $M_t$.
\item Extract $L_t$ and $\sqrt{\mathbf{S}_t}$ from $U$, partitioned as in Eq.\ \eqref{eq:sqrtq}:
\begin{align}
L_t & \gets U_{00}^{-1} U_{01}, & \sqrt{\mathbf{S}_t} & \gets U_{11},
\end{align}
where back-substitution can be used for the left-multiplication by $U_{00}^{-1}$.
\end{enumerate}

The backward pass is initialized by setting $\sqrt{\mathbf{S}_\ell}$. By Eq.\ \eqref{eq:sell}, $\mathbf{S}_\ell = M_\ell^T M_\ell$ with $M_\ell = \sqrt{W_\ell}\, R_\ell$, so $\sqrt{\mathbf{S}_\ell}$ follows directly from the QR-factorization of $M_\ell$.

The forward pass is exactly as in Section \ref{sec:forwardpass}, with the difference that the value $-d > 0$ can be obtained by summing the squared length of the rightmost column of $U_{01}$ over all time steps during the square root backward pass.\footnote{Evaluating the inequality of Eq.\ \eqref{eq:armijo} requires computing the difference $j(\alpha) - j^\star$ of two potentially large numbers, which can cause catastrophic cancellation. By recognizing that the total cost of a trajectory is of the form $j = \frac{1}{2}\|\mathbf{\mathbf{j}}\|^2$, where $\mathbf{j}$ is the vector that stacks $[\sqrt{W_t}\, \mathbf{r}_t(\mathbf{x}_t, \mathbf{u}_t)]$ for all $t$, this difference can be computed accurately as $j(\alpha) - j^\star = \frac{1}{2}(\mathbf{j}(\alpha) + \mathbf{j}^\star) \cdot (\mathbf{j}(\alpha) - \mathbf{j}^\star)$.}

\subsection{Computational complexity}
The relationship between standard and square root iLQR mirrors the classical tradeoff between the normal equations and QR-decomposition in linear least squares. A careful accounting of floating-point operations per backward-pass step\textemdash taking into account symmetry and triangularity\textemdash yields the following (omitting sub-cubic terms), where $k$, $m$, and $n$ are the dimensions of the residual, control, and state, respectively:
\begin{align}
\text{Square root iLQR:} & \quad \sim 2k(m + n)^2 + \tfrac{7}{3}n^3 + 3mn^2 + m^2n - \tfrac{2}{3}m^3, \\
\text{Standard iLQR:} & \quad \sim k(m + n)^2 + 3n^3 + 5mn^2 + 2m^2n + \tfrac{1}{3}m^3.
\end{align}
The leading $k$-dependent term doubles, so the overhead of square root iLQR is at most a factor of 2\textemdash the classic price of QR over normal equations \cite{golub}. The remaining terms, however, are smaller, and the two methods break even when
\begin{align}
k = \frac{m^3 + m^2 n + 2 m n^2 + \tfrac{2}{3} n^3}{(m+n)^2}. \label{eq:crossover}
\end{align}
For $m = n$, this crossover occurs at $k = \tfrac{7}{6}n$, i.e., close to the minimum residual dimension $k = m$ required by the full-rank condition on $\nabla_\mathbf{u}\mathbf{r}_t$. In practice the overhead is modest for problems where $k$ is on the order of $m + n$: on the bimanual manipulation problem of Section~\ref{sec:experiments} a complete trajectory optimization is within $3\%$ of the standard solver on a $3.0$\,GHz Xeon E5-2687W v4.

\section{Constrained iLQR}\label{sec:constraints}

In many practical applications, the trajectory must satisfy additional constraints beyond the dynamics in order to, for example, avoid collisions or respect actuator limits. An  effective approach to handling constraints in iLQR is the \emph{augmented Lagrangian} method~\cite{altro}, which converts the constrained problem into a sequence of unconstrained problems~\cite{bertsekas}. We extend the trajectory optimization problem of Eq.\ \eqref{eq:totalcost} by adding inequality constraints at each time step:
\begin{align}
\text{minimize} &\quad c_\ell(\mathbf{x}_\ell) + \sum_{t = 0}^{\ell - 1} c_t(\mathbf{x}_t, \mathbf{u}_t) \nonumber \\
\text{subject to} &\quad \mathbf{x}_{t+1} = \mathbf{f}_t(\mathbf{x}_t, \mathbf{u}_t),  \quad \mathbf{x}_0 = \mathbf{x}_0^\star, \quad \mathbf{g}_t(\mathbf{x}_t, \mathbf{u}_t) \leq \mathbf{0}, \label{eq:constrainedproblem}
\end{align}
where $\mathbf{g}_t(\mathbf{x}_t, \mathbf{u}_t) \in \mathbb{R}^p$ is a vector of $p$ inequality constraints at time step $t$.

The augmented Lagrangian approach converts constraints into cost terms by introducing \emph{Lagrange multipliers} $\boldsymbol{\lambda}_t \in \mathbb{R}^p$ and a \emph{penalty parameter} $\mu > 0$. At each stage $t$, we form the augmented cost $c_t^+$:
\begin{align}
c^+_t(\mathbf{x}_t, \mathbf{u}_t) & = \frac{\mu}{2} \left\|\max\{\mathbf{0}, \boldsymbol{\lambda}_t / \mu + \mathbf{g}_t(\mathbf{x}_t, \mathbf{u}_t)\}\right\|^2, \label{eq:alcost}
\end{align}
and this is simply \emph{added} to the original running cost $c_t$ at each time step. Equality constraints $\mathbf{h}_t(\mathbf{x}_t, \mathbf{u}_t) = \mathbf{0}$ can be handled analogously by omitting the maximization with zero (i.e., treating the constraints as always active).

The optimization proceeds as an outer--inner loop that repeats until the constraints are satisfied to within a desired tolerance:
\begin{enumerate}
\item \textbf{Inner loop.} Run iLQR to convergence on the unconstrained problem with the augmented cost $c_t + c^+_t$, treating $\boldsymbol{\lambda}_t$ and $\mu$ as constants.
\item \textbf{Outer loop.} Update the Lagrange multipliers and increase the penalty:
\begin{align}
\boldsymbol{\lambda}_t & \gets \max\{\mathbf{0}, \boldsymbol{\lambda}_t + \mu \mathbf{g}_t(\mathbf{x}_t, \mathbf{u}_t) \}, & \mu & \gets \phi\, \mu \quad (\phi > 1). \label{eq:lambdaupdate}
\end{align}
\end{enumerate}


The penalty term of Eq.~\eqref{eq:alcost} is a weighted square, and therefore the Gauss-Newton formulation applies directly. We linearize the constraint function $\mathbf{g}_t$ and form the cost Hessian from Jacobians alone, yielding the approximation:
\begin{align}
c^+_t(\mathbf{x}_t, \mathbf{u}_t) & \approx \frac{\mu}{2} \begin{bmatrix}\delta\mathbf{u}_t \\\delta\mathbf{x}_t \\ 1 \end{bmatrix}^T\! G_t^T D_t G_t\begin{bmatrix}\delta\mathbf{u}_t \\\delta\mathbf{x}_t \\ 1 \end{bmatrix}\!, & G_t & = \left. \begin{bmatrix} \nabla_\mathbf{u} \mathbf{g}_t & \nabla_\mathbf{x} \mathbf{g}_t & \boldsymbol{\lambda}_t/\mu\!+\!\mathbf{g}_t \end{bmatrix} \right\rvert_{\mathbf{x}^\star_t, \mathbf{u}^\star_t}\!,\nonumber
\end{align}
where $D_t = \diag(\boldsymbol{\lambda}_t / \mu + \mathbf{g}_t(\mathbf{x}_t^\star, \mathbf{u}_t^\star) > \mathbf{0}) \succeq 0$ is a diagonal  and idempotent matrix with a $0$ or $1$ on the diagonal for each constraint indicating whether it is active.
The resulting augmented cost matrix is positive-semidefinite by construction:
\begin{align}
\mathbf{C}_t & = R_t^T W_t R_t + \mu G_t^T D_t G_t \succeq 0.
\label{eq:augmented_cost}
\end{align}

The square root backward pass of Section~\ref{sec:sqrt} therefore applies unchanged but for one additional row block in $M_t$:\footnote{Note that $D_t$ is symmetric and idempotent. So, $D_t^T D_t = D_t^2 = D_t$.}
\begin{align}
\mathbf{Q}_t & = M_t^T M_t, & M_t & = \begin{bmatrix}
\sqrt{W_t}\, R_t \\
\sqrt{\mu}\,D_t G_t \\
\sqrt{\mathbf{S}_{t + 1}}\, F_t
\end{bmatrix}, \label{eq:matrixmconstrained}
\end{align}
and correspondingly $\smash{M_\ell = \bigl[\begin{smallmatrix} \sqrt{W_\ell} \, R_\ell \\ \sqrt{\mu}\, D_\ell G_\ell \end{smallmatrix}\bigr]}$ for the initialization.

The augmented Lagrangian framework is effective but introduces a tension with numerical precision. The larger the penalty parameter $\mu$, the faster the outer loop converges \cite{bertsekas}\textemdash but since $\mu$ enters as a weight in the cost Hessian, the eigenvalues of the cost-to-go matrix $\mathbf{S}_t$ scale with $\mu$, while eigenvalues associated with the unconstrained objectives remain at their original magnitudes. This results in a large condition number for the matrix $\mathbf{S}_t$, and may cause swamping or cancellation in the repeated matrix manipulations within the backward pass. 
In practice, rounding first erodes the accuracy of the computed policy and then, once $Q_\mathbf{uu}$ is indefinite or singular, prevents the standard backward pass from producing a policy at all. It is against exactly this that the square root formulation of Section~\ref{sec:sqrt} is of benefit, and Section~\ref{sec:experiments} measures by how much.

\section{Numerical Experiments}\label{sec:experiments}

We validate the practical benefits of square root iLQR on a constrained trajectory optimization problem within the augmented Lagrangian framework and its associated conditioning challenges.

\begin{figure*}[t]
\centering
\includegraphics[width=\textwidth]{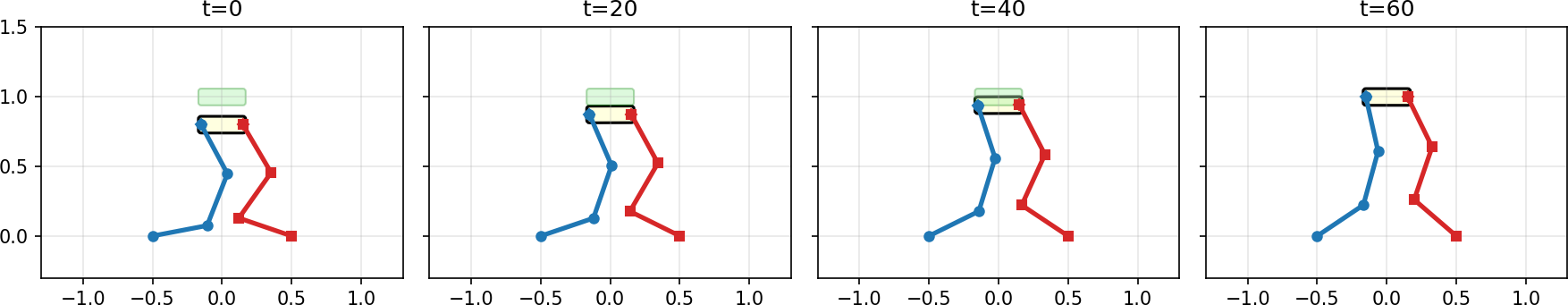}
\caption{Bimanual manipulation trajectory. Two 3-link arms transport a rigid object from the initial pose ($t\!=\!0$) to the goal at $(0, 1, 0)$ (green rectangle), shown at four keyframes. This trajectory was obtained using square root Gauss-Newton iLQR within an augmented Lagrangian loop.}
\label{fig:overlay}
\vspace{-3pt}
\end{figure*}

Consider two planar 3-link robot arms that must cooperatively transport a rigid object from an initial pose to a goal pose (Fig.~\ref{fig:overlay}). Each arm has link lengths $l = 0.4$\,m, with bases at $(\pm 0.5, 0)$, and grips the object at a body-frame offset of $(\pm 0.15, 0)$. The state $\mathbf{x} \in \mathbb{R}^6$ consists of three joint angles per arm ($n = 6$), and the control $\mathbf{u} \in \mathbb{R}^6$ consists of joint velocities. The dynamics are $\mathbf{x}_{t+1} = \mathbf{x}_t + \Delta t\, \mathbf{u}_t$ with $\Delta t = 0.05$\,s over $\ell = 60$ steps. The residual $\mathbf{r}_t$ comprises a goal-tracking term at the terminal step and velocity-penalizing terms at every step. Since the object pose is recovered from the two grasp points, closing the kinematic chain needs only a single scalar equality constraint $g(\mathbf{x}) = 0$, requiring the end-effectors to stay the correct distance apart, which we enforce via the augmented Lagrangian framework.

The key practical advantage of square root iLQR emerges when the penalty parameter $\mu$ grows large. To isolate the effect of numerical conditioning from other convergence factors, we first solve the problem to convergence at moderate $\mu$ using the full augmented Lagrangian loop, obtaining a well-converged trajectory. We then evaluate a \emph{single backward pass} at artificially elevated values of $\mu$, comparing the feedback gain matrices $L_t$ produced by each method against a reference.  Both methods execute an identical sequence of operations in exact arithmetic, so any difference in $L_t$ is pure numerical error. The value of $\mu$ here serves as a proxy for the condition numbers encountered among the cost-to-go matrices, which grow linearly in $\mu$.

\begin{figure}[t]
\centering
\includegraphics[width=\columnwidth]{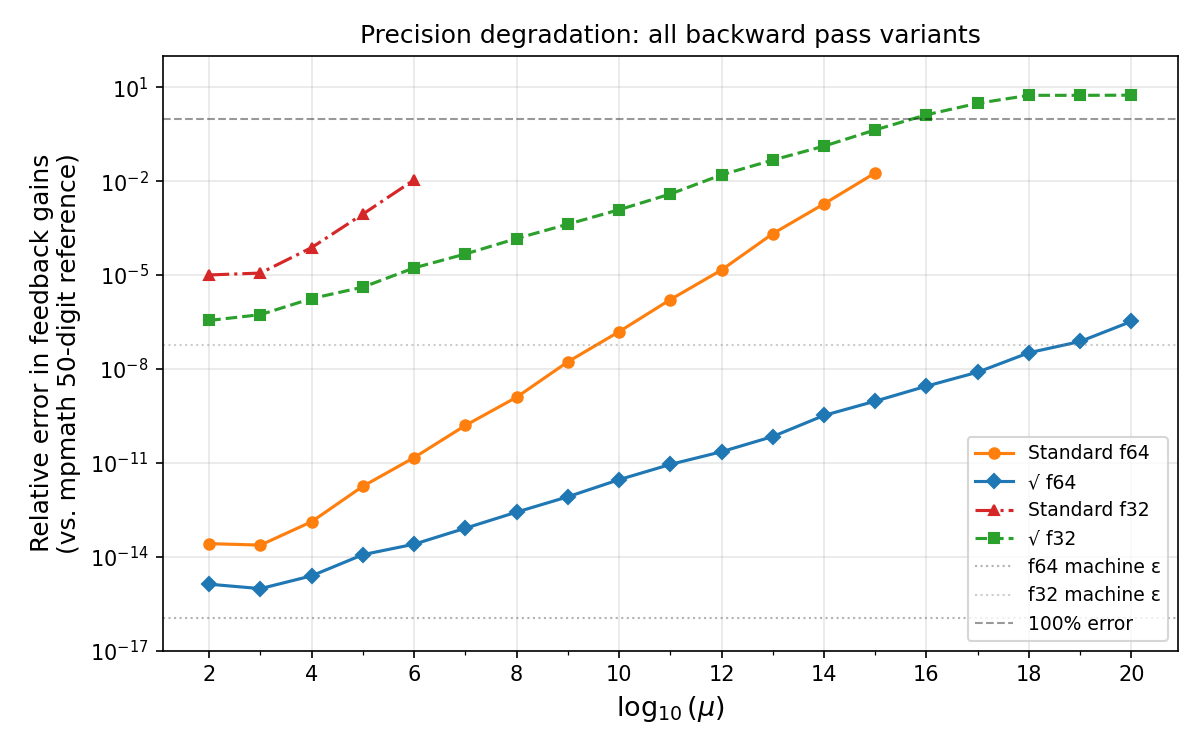}
\caption{Maximum relative error $\max_t \lVert L_t - L_t^{\star}\rVert_F / \lVert L_t^{\star}\rVert_F$ in feedback gains for four backward pass variants, measured against a 50 decimal digit arbitrary-precision reference $L_t^{\star}$ computed with Python's \texttt{mpmath} library: the standard backward pass in f64 (orange), the square root backward pass in f64 (blue), the standard backward pass in f32 (red dash-dot), and the square root backward pass in f32 (green dashed).}
\label{fig:precision}
\vspace{-3pt}
\end{figure}

Fig.~\ref{fig:precision} shows the relative error in feedback gains as a function of $\mu$ for four backward pass variants, all measured against a reference computed at 50 decimal digits of precision using Python's \texttt{mpmath} library. All variants degrade with striking regularity, with the standard iLQR variants losing approximately one significant digit per decade of $\mu$, and the square root variants losing significant digits at half that rate. This is consistent with the fact that Cholesky factors have condition numbers that are the square root of those of the full matrices.

For the standard f64 variant (orange), barely two significant digits survive at $\mu = 10^{15}$, and at $\mu = 10^{16}$ the Cholesky decomposition of $Q_\mathbf{uu}$ fails outright. The square root f64 variant (blue) starts an order of magnitude closer to machine epsilon and still retains six to seven significant digits at $\mu = 10^{20}$. Repeating the experiment in f32 arithmetic confirms that the pattern is intrinsic to the formulation and not an artifact of f64: the two f32 curves reproduce the f64 pair, slopes and low-$\mu$ separation alike, translated upward in error by roughly the eight decades that separate the two machine epsilons. The standard f32 pass (red dash-dot) consequently fails already at $\mu = 10^{7}$.

Whereas the standard variants fail for too large values of $\mu$ (the matrix $Q_\mathbf{uu}$ becomes indefinite or singular), the square root variants degrade silently. They produce meaningless results when $\mu$ grows too large, as can be seen from the (green dashed) curve of the square root f32 variant for $\mu = 10^{15}$ and upwards, without any indication of error. This is because numerical positivity is guaranteed when using Cholesky factors (although they may still become singular). A health check comes free, however, since $U_{00} = \sqrt{Q_\mathbf{uu}}$ is already triangular. The squared ratio of its largest to its smallest diagonal magnitude bounds the condition number of $Q_\mathbf{uu}$ from below. Its smallest diagonal entry tends to zero exactly as $\nabla_\mathbf{u}\mathbf{r}_t$ loses column rank. We flag a step once that entry falls below $\sqrt{\epsilon}$ times the largest, where $\epsilon$ is the machine epsilon; by that point half the available digits are gone. A flagged step can trigger regularization or early termination.

The advantage of square root iLQR does not depend on problem size. A second problem supplies a size knob: a planar chain of $N$ links unfurling toward a goal outside its circular workspace, with an inequality on each constrained tip that activates as that tip reaches the boundary. We swept three parameters, doubling each in turn at fixed arm length: the state and control dimension ($N = 4$ to $32$, four tips constrained), the horizon ($\ell = 20$ to $160$, $N = 16$) and the number of constrained tips ($4$ to $32$, on a $32$-link chain). Absolute failure thresholds do move, by up to two decades over the dimension range. Both backward passes move together, however, so the margin between them survives: in every configuration the square root formulation buys five to six extra decades of $\mu$ in f32.

A single backward pass isolates round-off, but whether the added precision changes the outcome of a complete solve is a separate question. We ran the full augmented Lagrangian loop from $24$ randomized initial trajectories with an f32 backward pass, tightening the constraint tolerance until each method failed. Both reach $10^{-6}$; at $10^{-7}$ the square root solver still succeeds from all $24$ starts against $3$, and at $10^{-8}$ from all $24$ against none, the successful runs agreeing on the same trajectory and final cost. In f64 the same divergence appears near $10^{-11}$. What the formulation buys is therefore not a better optimum but how much further the outer loop can be driven before the backward pass gives out, allowing as a result for much tighter constraint tolerances.

\section{Discussion and Conclusion}
Many trajectory optimization problems can be formulated with weighted sums-of-squares cost functions in practice, including those with (in)equality constraints as we have shown in this paper. For such problems, we have presented an elegant and easy-to-implement square root formulation of Gauss-Newton iLQR that roughly doubles the effective numerical precision at negligible additional cost\textemdash under $3\%$ of computation time over a complete solve. Our approach reduces each step of the backward pass to computing a single QR-decomposition of a fixed-size matrix, and as such qualitatively improves upon previous square root iLQR variants \cite{geoffroy,altro}. The authors believe that when it is applicable, there is no reason not to use the square root Gauss-Newton iLQR over standard Gauss-Newton iLQR, analogous to the fact that there are (almost) no downsides to using a square root Kalman filter over a standard one in any production system.

The relative success of square root variants of the Kalman filter compared to square root variants of iLQR may be explained by the fact that when posed as an optimization problem, any linear-Gaussian estimation problem also adheres to the Gauss-Newton cost structure, and can therefore benefit from equally elegant square root solutions. In contrast, iLQR and related algorithms are typically presented in the context of the most general problem definition, which does not readily lend itself to a square root implementation. 

Even so, Appendix~\ref{sec:general} sketches a possible square root formulation for general cost functions, though such a formulation is necessarily less elegant. Ultimately, it is the Gauss-Newton cost structure that makes the square root formulation of iLQR that we introduced in this paper so pleasingly simple, and we hope that the success square root variants have enjoyed in the estimation space can be extended to the area of trajectory optimization.

\appendix
\section{Square Root iLQR for General Cost Functions} \label{sec:general}
For completeness, we sketch how a square root formulation of iLQR can be constructed for general cost functions as well. In this case $\sqrt{S_\mathbf{xx}}$, $\mathbf{s}_\mathbf{x}$ and $s$ (see Eq.\ \eqref{eq:costogo}) are propagated separately. For a square root formulation to work, a (possibly non-triangular) square root of the Hessian block of the cost matrix $\mathbf{C}_t$ must be available in each step. Since this Hessian is in general indefinite, such a square root does not exist. Therefore, we propose to use the square root of the nearest positive-semidefinite approximation of this Hessian. That is, we eigendecompose the Hessian of the cost function:
\begin{align}
\left. \begin{bmatrix} \nabla^2_\mathbf{uu} c_t  & \nabla^2_\mathbf{ux} c_t \\ \nabla^2_\mathbf{xu}c_t & \nabla^2_\mathbf{xx} c_t\end{bmatrix} \right \rvert_{\mathbf{x}^\star_t,\mathbf{u}^\star_t} = P \Lambda P^T,
\end{align}
where $P$ is orthogonal and $\Lambda$ is a diagonal matrix of eigenvalues, then set negative eigenvalues to zero and positive ones to their square root, and let:
\begin{align}
\tilde{M}_t & = \begin{bmatrix} \sqrt{\max\{\Lambda, 0\}}\, P^T \\
\sqrt{S_\mathbf{xx}}_{,t + 1} \begin{bmatrix} \nabla_\mathbf{u} \mathbf{f}_t & \nabla_\mathbf{x} \mathbf{f}_t \end{bmatrix}
\end{bmatrix},
\end{align}
such that we have
$\bigl[\begin{smallmatrix} Q_\mathbf{uu} & Q_\mathbf{ux} \\ Q_\mathbf{ux}^T & Q_\mathbf{xx}\end{smallmatrix}\bigr] = \tilde{M}_t^T \tilde{M}_t$. If we let $L_t = \begin{bmatrix} K_t & \mathbf{k}_t \end{bmatrix}$, we have for the Cholesky factor $\tilde{U}$ of the Hessian block of $\mathbf{Q}_t$:
\begin{align}
\tilde{U} = \begin{bmatrix} \tilde{U}_{00} & \tilde{U}_{01} \\ 0 & \tilde{U}_{11} \end{bmatrix} = \sqrt{\begin{bmatrix} Q_\mathbf{uu} & Q_\mathbf{ux} \\ Q_\mathbf{ux}^T & Q_\mathbf{xx}\end{bmatrix}} = \begin{bmatrix} \sqrt{Q_\mathbf{uu}} & \sqrt{Q_\mathbf{uu}}\, K_t \\ 0 & \sqrt{S_\mathbf{xx}}_{,t} \end{bmatrix}, \label{eq:hessiancholesky}
\end{align}
and it is obtained as the upper-triangular factor of the QR-factorization of $\tilde{M}_t$. From $\tilde{U}$, we can extract $K_t = \tilde{U}_{00}^{-1}\tilde{U}_{01}$ and $\sqrt{S_\mathbf{xx}}_{,t} = \tilde{U}_{11}$. Further, we have:
\begin{align}
\mathbf{q}_\mathbf{u} & = \nabla_\mathbf{u}^T \mathbf{f}_t\, \mathbf{s}_{\mathbf{x}, t+1} + \nabla_\mathbf{u}c_t, &
\mathbf{q}_\mathbf{x} & = \nabla_\mathbf{x}^T \mathbf{f}_t\, \mathbf{s}_{\mathbf{x}, t+1} + \nabla_\mathbf{x}c_t, &
q & = s_{t + 1} + 2c_t,
\end{align}
from which the other desired quantities\textemdash the feedforward term $\mathbf{k}_t$ and the gradient $\mathbf{s}_\mathbf{x}$ and value $s$ of the updated cost-to-go function\textemdash can be computed:
\begin{align}
\tilde{\mathbf{q}}_\mathbf{u} & = \tilde{U}_{00}^{-T} \mathbf{q}_\mathbf{u}, & \mathbf{k}_t & = \tilde{U}_{00}^{-1} \tilde{\mathbf{q}}_\mathbf{u}, & \mathbf{s}_{\mathbf{x}, t} & = \mathbf{q}_\mathbf{x} - \tilde{U}_{01}^T \tilde{\mathbf{q}}_\mathbf{u}, & s_t & = q - \tilde{\mathbf{q}}_\mathbf{u}^T \tilde{\mathbf{q}}_\mathbf{u}.
\end{align}
The left-multiplications by $\tilde{U}_{00}^{-T}$ and $\tilde{U}_{00}^{-1}$ can be performed by forward- and back-substitution, respectively.

The above construction is perfectly workable, but it reintroduces both squared quantities and regularization logic\textemdash precisely what the Gauss-Newton structure of the main text dispenses with.



\newpage


\end{document}